\documentclass[12pt]{article}
\usepackage{amssymb}
\usepackage{amsfonts}
\usepackage{color}
\usepackage[pdftex]{graphicx}

\begin{document}
\begin{center}
\Large
Boole's formula as a consequence of Lagrange's Interpolating Polynomial theorem
\end{center}

\bigskip

\begin{center}
{\bf{Cosmin Pohoata}}

{\it{"Tudor Vianu" National College of Informatics, Bucharest, Romania}}

\verb+pohoata_cosmin2000@yahoo.com+
\end{center}

\bigskip

\bigskip

\begin{center}
{\bf{Abstract}}
\end{center}
We present a slightly more general version of Boole's additive formula for factorials as a simple consequence of Lagrange's Interpolating Polynomial theorem.

\bigskip

\bigskip 

\bigskip

{\bf{1. Introduction}} 

\bigskip

\bigskip

In the first chapter of [3], Boole defines, for all real valued function of one real variable $f(x)$, the {\it{first difference}} of $f(x)$ (with respect to the increment 1) as $\Delta f(x) = f(x+1)-f(x)$. He then defines, for all integers $n \geq 2$, the $n$-th difference by the recursive formula $\Delta^{n}f(x) = \Delta\Delta^{n-1}f(x)$. This enables him to prove by induction (see [3], p. 5, (2)) that, for all positive integers $n$,
$$\Delta^{n}x^{n} = n!.$$
Later on, with the help of this formula, he derives the following identity, known nowadays as the {\it{Boole additive formula for factorials}} (see [3], p. 20, (6)),
$$\sum_{k=1}^{n}{(-1)^{n-k}{{n}\choose k}k^{n}}=n!$$

In [1] R. Anglani and M. Barile present two proofs of this identity, one completely analytical, and the other by making use of an ingenuous combinatorial argument. Here, we will give an immediate proof of a more general version of this identity by making use of Lagrange's Interpolating Polynomial theorem (see, for example, [2]).

\bigskip

\bigskip

\bigskip

{\bf{2. Main result}}

\bigskip

\bigskip

{\bf{Proposition}}. Consider $p(x)=a_{0}x^{n}+a_{1}x^{n-1}+\ldots a_{n-1}x+a_{n}$, an arbitrary polynomial of degree $n$ with real coefficients. For any real numbers $a$, $b$, with $b \neq 0$,
$$\sum_{k=0}^{n}{(-1)^{n-k}{{n}\choose k}p(a+kb)}=a_{0} \cdot b^{n} \cdot n!.$$

\bigskip

{\it{Proof}}. Since $p$ has degree at most $n$, according to the Lagrange Interpolating Polynomial theorem,
$$p(x)=\sum_{k=0}^{n}{p(a+kb)\prod_{0 \leq j\neq k \leq n}{\frac{x-a-jb}{(k-j)b}}}.$$
By identifying the leading coeficients on both sides of the above equality,
$$a_{0}=\sum_{k=0}^{n}{p(a+kb)\prod_{0 \leq j\neq k \leq n}{\frac{1}{(k-j)b}}}=\frac{1}{b^{n}}\sum_{k=0}^{n}{\frac{(-1)^{n-k}p(a+kb)}{k!(n-k)!}}$$
$$=\frac{1}{n!b^{n}}\sum_{k=0}^{n}(-1)^{n-k}{{n}\choose k}p(a+kb). \ \ \ \ \ \ \ \ \ \ \ \ \ \ \ \ \ \ \ \ \ \ \ \ \ \ \ \ \ \ \ $$
Therefore,
$$\sum_{k=0}^{n}{(-1)^{n-k}{{n}\choose k}p(a+kb)}=a_{0} \cdot b^{n} \cdot n!$$
$\hfill \square$

\bigskip

Obviously, for $a=0$, $b=1$, and $p(x)=x^{n}$, Proposition 1 is equivalent to Boole's formula. Moreover, it yields the well-known identity
$$\sum_{k=1}^{n}{(-1)^{n-k}{{n}\choose k}k^{m}}=0,$$
which holds for $m=0,\ 1,\ \ldots,\ n-1$.

\bigskip

\bigskip

\bigskip

{\bf{References}}

\bigskip

1. R. Anglani, M. Barile, Two very short proofs of a combinatorial identity, The Electronic Journal of Combinatorial Number Theory 5 (2005).

\bigskip

2. Archer, Branden and Weisstein, Eric W. "Lagrange Interpolating Polynomial", available at

\verb+http://mathworld.wolfram.com/LagrangeInterpolatingPolynomial.html+.

\bigskip

3. G. Boole, Calculus of Finite Differences, (J. F. Moulton, ed.), 4th Edition. Chelsea,
New York. (n.d.)

\bigskip

\end{document}